\title{\bf{Some infinite sequences of canonical covers of degree 2}}
\author{
	NGUYEN BIN\\
}
\date{\today}

\newcommand{\Addresses}{{
		\bigskip
		\footnotesize
			\text{Center for Mathematical Analysis, Geometry and Dynamical Systems}\par\nopagebreak
			\text{Departamento de Matem\'{a}tica}\par\nopagebreak
			\text{Instituto Superior T\'{e}cnico}\par\nopagebreak
			\text{Universidade de Lisboa}\par\nopagebreak
			\text{Av. Rovisco Pais}\par\nopagebreak
			\text{1049-001 Lisboa}\par\nopagebreak
			\text{Portugal.}\par\nopagebreak
		\textit{E-mail address}: \texttt{nguyenbin@tecnico.ulisboa.pt}
				
	}}

\newcommand\blfootnote[1]{%
	\begingroup
	\renewcommand\thefootnote{}\footnote{#1}%
	\addtocounter{footnote}{-1}%
	\endgroup
}

\date{\today}

\documentclass[10pt]{article}
\usepackage{amsmath, amsthm, amssymb, mathrsfs, amscd, amsthm, amscd,amsfonts,calligra,mathrsfs,adjustbox,lipsum}

\usepackage{indentfirst,url,xypic}
\usepackage[all]{xy}
\usepackage{url}
\usepackage{tikz}
\usepackage[colorlinks,plainpages]{hyperref}
\hypersetup{
	colorlinks=true,
	linkcolor=blue,
	filecolor=magenta,      
	urlcolor=cyan,
}

\DeclareMathOperator{\degree}{deg}
\DeclareMathOperator{\image}{im}
\DeclareMathOperator{\Picard}{Pic}

\newtheorem{Theorem}{Theorem}
\newtheorem{Proposition}{Proposition }

\theoremstyle{remark}
\newtheorem{Remark}{Remark}

\newtheorem*{Acknowledgments}{ACKNOWLEDGMENTS}

\newcommand{\MSC}{\textbf{Mathematics Subject Classification (2010):}}

\newcommand{\Key}{\textbf{Key words:}}

\begin{document}
\maketitle
\begin{abstract}  In this note, we construct three new infinite families of surfaces of general type with canonical map of degree 2 onto a surface of general type. For one of these families the canonical system has base points. 
\end{abstract}

\blfootnote{\MSC{ 14J29}.}
\blfootnote{\Key{ Surfaces of general type, Canonical maps, Abelian covers.}}

\section{Introduction}
 Let $ X $ be a minimal smooth complex surface of general type and denote by $ \xymatrix{\varphi_{\left| K_X\right| }:X \ar@{.>}[r] & \mathbb{P}^{p_g\left( X\right)-1}} $ the canonical map of $ X $, where $ p_g\left( X\right) = \dim\left( H^0\left( X, K_{X}\right) \right) $ is the geometric genus and $ \left| K_{X} \right| $ is the canonical system of $ X $. The study of the degree of the canonical map for surfaces of general type has a long history and was inspired by the following well-known work of A. Beauville \cite{MR553705}.
\begin{Theorem} \cite[ \rm Theorem 3.1]{MR553705} \label{Theorem of Beaville}
	If the image $ \Sigma:=\image\left( \varphi_{\left | K_X\right| }\right)  $ is a surface, then either 
	\begin{enumerate}
		\item $ p_g\left( \Sigma\right) =0, $ or
		\item $ \Sigma $ is a canonical surface (in particular, $ p_g\left( \Sigma\right) =p_g\left( X \right) $).
	\end{enumerate}
	Moreover, in case $ \left( 1\right),  $ $ d:= \degree\left( \varphi_{\left | K_X\right| }\right) \le 36, $ and in case $ \left( 2\right) $, $ d \le 9 $.
\end{Theorem}
In the same paper \cite{MR553705}, A. Beauville also proved that in the second case if $\chi(\mathcal O_X)\geq 14$ the degree of the canonical map is less than or equal to $3 $. \\

Only a few  families of surfaces with arbitrarily high geometric genus in the second case are known. The first construction due to A. Beauville (see \cite{MR927956} or \cite{MR1659365}) gives an infinite series of surfaces with $ d = 2 $ and irregularity $ q=2 $. Later, M.M. Lopes and R. Pardini proved that these are almost the only examples satisfying $ K^2 = 6p_g - 14$ and $ q \ge 2 $ (see \cite{MR1659365}). By means of a construction analogous to that of A. Beauville's, C. Ciliberto, R. Pardini and F. Tovena created two more series of surfaces: one with $ d = 2 $ and $ q=2 $, and another with $ d = 2 $ and $ q=3 $ (see \cite{MR1822412}). A few years later they found three more unlimited families of surfaces with $ d = 2$ and $ q=0 $ (see \cite{MR1960801}). These six  families of surfaces all satisfy $ K^2 < 7p_g $. Furthermore, these surfaces all possess canonical linear systems which are base point free.\\

 As far as we know, all the known examples of the second case of Theorem \ref{Theorem of Beaville} have base point free canonical system. In this note, we construct three new families of surfaces satisfying the second case of Theorem \ref{Theorem of Beaville} and for one of the families the canonical system has base points. These three infinite families of surfaces all satisfy either $ K^2=8p_g-24  $ or $ K^2=8p_g-18  $. The following theorem is the main result of this note:
\begin{Theorem}\label{the main theorem}
	Let $ m, n $ be integer numbers such that $ m,n \ge 2 $. Then there exist minimal surfaces of general type $ X $ whose canonical map is of degree 2 onto a surface of general type $ \Sigma $ with the following invariants
	$$
	\begin{tabular}{|c| c| c| c| c|}
	\hline
	$ K_X^2 $ &$ p_g\left( X\right) $ &$ q\left( X\right) $ & $ \degree\left( \Sigma\right)  $& $\left| K_X \right| $ is base point free \\
	\hline
	$ 16mn $&$ 2mn + 3 $&$ 0 $&$ 8mn $& yes \\
	\hline
	$ 16mn-8 $&$ 2mn + 2 $&$ 0 $&$ 8mn -4 $& yes \\
	\hline
	$ 16mn-2 $&$ 2mn + 2 $&$ 0 $&$ 8mn -2 $& no \\
	\hline
	\end{tabular} 
	$$
	
\end{Theorem}

The idea of our constructions is the following. \\
We construct the surfaces in the first row of Theorem \ref{the main theorem} by taking a specific $ \mathbb{Z}_{2}^3- $cover $ X $ of $ \mathbb{P}^1 \times \mathbb{P}^1 $. In particular, we take the cover to be in such a way that one of the $ \mathbb{Z}_{2}- $quotients $ Z $ of $X$  is a surface of general type whose only singularities are nodes (i.e. points of type $A_1$). Moreover, we require that: 
\begin{enumerate}
	\item $ p_g\left( Z\right) = p_g\left( X\right)  $,
	\item the canonical map $ \varphi_{\left| K_{Z} \right| } $ is birational.	
\end{enumerate}
\noindent
To obtain  the other surfaces described in Theorem \ref{the main theorem}, we modify this $ \mathbb{Z}_{2}^3- $cover by imposing an ordinary quadruple point on the branch locus in such a way that the two conditions above hold.

\section{$ \mathbb{Z}_{2}^3- $coverings}
The construction of abelian covers was studied by R. Pardini in \cite{MR1103912}. 
For details about the building data of abelian covers  we refer the reader to Section 1 and Section 2 of R. Pardini's work (\cite{MR1103912}). 

We will denote by  $ \chi_{j_1j_2j_3} $ the character of $ \mathbb{Z}_{2}^3 $ defined by
\begin{align*}
	\chi_{j_1j_2j_3}\left( a_1,a_2,a_3\right): =  e^{\left( \pi a_1j_1\right) i}e^{\left( \pi a_2j_2\right) i}e^{\left( \pi a_3j_3\right) i}
\end{align*}
for all $ j_1,j_2,j_3,a_1,a_2,a_3\in \mathbb{Z}_2 $. 

From \cite[\rm Theorem 2.1]{MR1103912} we can define $ \mathbb{Z}_{2}^3- $covers as follows:
\begin{Proposition} \label{Construction of cover}
	Given $ Y $ a smooth projective surface, let $ L_{\chi} $ be divisors of $ Y $ such that $ L_{\chi} \not\equiv \mathcal{O}_Y $ for all nontrivial characters $ \chi $ of $ \mathbb{Z}_{2}^3  $ and let $ D_{\sigma} $ be effective divisors of  $ Y $ for all $ \sigma \in \mathbb{Z}_{2}^3 \setminus \left\lbrace \left(0,0,0 \right)  \right\rbrace  $ such that the branch divisor $ B:=\sum\limits_{\sigma \ne 0}{D_{\sigma}} $ is reduced. Then $ \left\lbrace L_{\chi}, D_{\sigma} \right\rbrace_{\chi,\sigma}$ is the building data of a $ \mathbb{Z}_{2}^3-$cover $ \xymatrix{f:X \ar[r]& Y} $ if and only if
	$$
	\begin{adjustbox}{max width=\textwidth}
	\begin{tabular}{l l l l l l l l }
	$ 2L_{100} $&$ \equiv $&$  $&$ $&$ D_{100} $&$ +D_{101 } $&$ +D_{110} $&$ +D_{111} $ \\
	$ 2L_{010} $&$ \equiv $&$ D_{010} $&$ +D_{011} $&$  $&$ $&$ +D_{110} $&$ +D_{111} $ \\
	$ 2L_{001} $&$ \equiv D_{001 } $&$ $&$ +D_{011} $&$$&$ +D_{101} $&$  $&$ +D_{111 } $ \\
	$ 2L_{110} $&$ \equiv $&$ D_{010 } $&$ +D_{011} $&$ +D_{100} $&$ +D_{101} $&$  $&$  $\\
	$ 2L_{101} $&$ \equiv D_{001} $&$ $&$ +D_{011} $&$ +D_{100} $&$  $&$ +D_{110} $&$  $ \\
	$ 2L_{011} $&$ \equiv D_{001} $&$ +D_{010} $&$  $&$  $&$ +D_{101} $&$ +D_{110} $&$  $ \\
	$ 2L_{111} $&$ \equiv D_{001} $&$ +D_{010} $&$  $&$ +D_{100} $&$ $&$ $&$ +D_{111 } $.
\end{tabular}
\end{adjustbox}
$$	
\end{Proposition}

By \cite[\rm Theorem 3.1]{MR1103912} if each $D_\sigma$ is smooth and $B $ is a simple normal crossings divisor, then the surface $X$ is smooth. \\

Also from \cite[\rm Lemma 4.2, Proposition 4.2]{MR1103912} we have:
\begin{Proposition}\label{invariants}
  If \hskip 2pt $ Y $ is a smooth surface and $ \xymatrix{f: X \ar[r]& Y} $ is a smooth $  \mathbb{Z}_{2}^3- $cover with building data $ \left\lbrace L_{\chi}, D_{\sigma} \right\rbrace_{\chi,\sigma}$, the surface $ X $ satisfies the following:
\begin{align*}
2K_X & \equiv f^*\left( 2K_Y + \sum\limits_{\sigma \ne 0} {D_{\sigma} } \right) \\
K^2_X &= 2\left( 2K_Y + \sum\limits_{\sigma \ne 0} {D_{\sigma} } \right)^2 \\
p_g\left( X \right) &=p_g\left( Y \right) +\sum\limits_{\chi \ne  \chi_{000}  }{h^0\left( L_{\chi} + K_Y \right)} \\
\chi\left( \mathcal{O}_X \right) &= 8\chi\left( \mathcal{O}_Y \right)  +\sum\limits_{\chi \ne \chi_{000}  }{\frac{1}{2}L_{\chi}\left( L_{\chi}+K_Y\right)}. 
\end{align*}
\end{Proposition}

\section{Construction}
\subsection{Construction and computation of invariants} \label{the original construction of surface with d = 2}
We will denote by $ F = \left\lbrace 0\right\rbrace \times \mathbb{P}^1  $ and $ G =   \mathbb{P}^1 \times \left\lbrace 0\right\rbrace$ the generators of $ \Picard\left( \mathbb{P}^1 \times \mathbb{P}^1  \right)  $. Let $ D_{100}$, $ D_{101} $ $\in \left|  2F +2G\right|  $, $ D_{110} \in \left|  2mF \right|  $ (with $m\geq 2$) and $ D_{111}$ $\in \left|  2nG \right|  $ (with $n\geq 2$) be smooth divisors of $  \mathbb{P}^1 \times \mathbb{P}^1 $ intersecting transversally  and such that no more than two of these divisors go through the same point. 

Let $ \xymatrix{f: X \ar[r] & \mathbb{P}^1 \times \mathbb{P}^1}  $ be the $ \mathbb{Z}^3_2- $cover with branch locus $ B = \sum\limits_{\sigma \ne 0}{D_{\sigma}} $,  where $ D_{\sigma} = 0 $ for the other  $ D_{\sigma}$. The building data is as follows:
$$	
\begin{tabular}{l r r}
	$ L_{100} \equiv$ & $ \left( m+2 \right)F   $ & $ +\left( n+2\right)G $\\
	$ L_{010} \equiv$ & $mF   $ &$ +nG$\\
	$ L_{110} \equiv$ & $2F   $ &$ +2G$\\
	$ L_{001} \equiv$ & $ F   $ & $ +\left( n+1\right)G $\\
	$ L_{101} \equiv$ & $ \left( m+1\right)F   $ & $ +G $\\
	$ L_{011} \equiv$ & $ \left( m+1\right)F   $ & $ +G $\\
	$ L_{111} \equiv$ & $ F   $ & $ +\left( n+1\right)G $.
\end{tabular} 
$$

The surface $X$ is smooth and satisfies
\begin{align*}
2K_X &\equiv f^*\left( 2mF + 2nG  \right). 
\end{align*}
\noindent
This implies that the surface $ X $ is minimal and of general type.  Furthermore, from Proposition \ref{invariants}, the surface $X$ has the following invariants:
\begin{align*}
K_X^2= 16mn, p_g\left( X\right) = 2mn+3, \chi\left( \mathcal{O}_X\right) =2mn+4, q\left( X\right)  = 0.   
\end{align*}

\subsection{The canonical map and the canonical image}\label{The canonical map and the canonical image}
\noindent
In this section, we show that the degree $d$ of the canonical map of $X$ is $ 2 $ and that the canonical image $ \Sigma $ is a surface of general type. We have the $ \mathbb{Z}_2^3- $equivariant decomposition 
\begin{align*}
H^{0}\left( X, K_X\right) = H^{0}\left( \mathbb{P}^1 \times \mathbb{P}^1, K_{\mathbb{P}^1 \times \mathbb{P}^1}\right) \oplus \bigoplus_{\chi \ne  \chi_{000}}{H^{0}\left( \mathbb{P}^1 \times \mathbb{P}^1, K_{\mathbb{P}^1 \times \mathbb{P}^1} +L_{\chi}\right)} 
\end{align*}
\noindent
where the group $ \mathbb{Z}_2^3 $ acts on $ H^{0}\left( \mathbb{P}^1 \times \mathbb{P}^1, K_{\mathbb{P}^1 \times \mathbb{P}^1} +L_{\chi}\right) $ via the character $ \chi $ (see \cite[ \rm Proposition 4.1c]{MR1103912}).

We consider the cyclic subgroup $ \Gamma:= \left\langle \left( 0,0,1\right) \right\rangle  $ of $ \mathbb{Z}_2^3 $. Let $ \Gamma^\perp $ denote the kernel of the restriction map $ \xymatrix{\left( \mathbb{Z}_2^3\right)^{*} \ar[r]&\Gamma^{*}} $, where $ \Gamma^{*} $ is the character group of $ \Gamma $. We have $ \Gamma^{\perp} = \left\langle \chi_{100},\chi_{010} \right\rangle  $. Because 
\begin{align*}
h^0\left( L_{ \chi} + K_X \right)  = 0
\end{align*}
for all $ \chi \notin \Gamma^\perp $, the subgroup $ \Gamma $ acts trivially on $ H^{0}\left( X, K_X\right) $. So the canonical map $ \varphi_{\left| K_X \right| } $ is the composition of the quotient map $ \xymatrix{X \ar[r]& Z:= X/\left\langle 0,0,1\right\rangle } $ with the canonical map $\varphi_{\left| K_{Z} \right| }$ of $ Z $ (see e.g. \cite[\rm Example 2.1]{MR1103913}). Thus, the following diagram commutes
$$
\xymatrix{X \ar[0,3]^{\mathbb{Z}_2^3}_f \ar[1,1] \ar[2,1]_{\varphi_{\left| K_X \right| }}&&& \mathbb{P}^1 \times \mathbb{P}^1\\
	&Z \ar[-1,2]_{\mathbb{Z}_2^2} \ar[1,0]^{\varphi_{\left| K_{Z} \right| }}&&\\
	&\mathbb{P}^{2mn+2}&&} 	
$$
\noindent
The intermediate surface $ Z $ is singular but since the singularities are ordinary double points, the formulas for the smooth case hold in this case. Since the intermediate surface $ Z $ is the fiber product of the double covers of $ \mathbb{P}^1 \times \mathbb{P}^1 $ branched on $ D_{100} + D_{101} $ and $ D_{110} + D_{111} $, respectively, the canonical class $ K_Z $ is the pullback of $ mF+nG $. In addition, the eigenspace decomposition of $ H^{0}\left( K_Z\right)  $ is as follows:
\begin{align*}
H^{0}\left( K_Z\right) = H^{0}\left( \mathcal{O}_{\mathbb{P}^1 \times \mathbb{P}^1}\right) \oplus H^{0}\left( \left( m-2\right) F + \left( n-2\right) G\right) \oplus H^{0}\left( mF + nG\right).
\end{align*}

\noindent
Because the linear system $ \left| K_Z\right|  $ contains the pullback of the very ample system $ \left| mF+nG \right|  $, there is a factorization $ \xymatrix{Z \ar[r]^{\varphi_{\left| K_{Z} \right| }} & Z^{'} \ar[r] & \mathbb{P}^1 \times \mathbb{P}^1} $. So if the map $ \varphi_{\left| K_{Z} \right| } $ is not birational, then the surface $ Z^{'} $ is birational to $ \mathbb{P}^1 \times \mathbb{P}^1 $ or to one of the three intermediate quotients $ Z/\left\langle g\right\rangle $, where $ g \in \mathbb{Z}_2^2 $ is a nonzero element. This does not happen since $ H^{0}\left( K_Z\right) $ has three nonzero eigenspaces. Thus, the canonical map $\varphi_{\left| K_{Z} \right| }$ is birational. Therefore, the canonical map of $X$ has degree $ 2 $ and the canonical image $ \Sigma $ is a surface of general type of degree $ 8mn $.

\subsection{Variations}
   We now consider variations of the construction of Section \ref{the original construction of surface with d = 2},  to obtain the surfaces listed in the second and third row of the table in Theorem \ref{the main theorem}. In order to achieve this, we impose an ordinary quadruple point on the branch locus and we resolve the singularities.
   
   \subsubsection{Variation 1}
   The first variation of the construction in Section \ref{the original construction of surface with d = 2} consists of allowing all the non zero $D_\sigma$ to pass through the same point.\\ 
      
   Namely we take a branch locus as before, except that now we require that the four divisors  $ D_{100}$, $ D_{101} $ $\in \left|  2F +2G\right|  $, $ D_{110} \in \left|  2mF \right|  $ (with $m\geq 2$) and $ D_{111}$ $\in \left|  2nG \right|  $ (with $n\geq 2$)  go through the same point $P$ in such a way that the branch locus $B$ has an ordinary quadruple point at $P$.  As before we require that   $ D_{100}$, $ D_{101} $, $ D_{110}$ and $D_{100}$ are smooth, that all intersections are transversal and that $P$ is the only point where more than two of the $D_\sigma$ meet. 	

   Denote by $ \overline{Y} $ the blow up of $ \mathbb{P}^1 \times \mathbb{P}^1 $ at $ P $, by $ \overline{F} $, $ \overline{G} $ the pullback of general fibers $ F $, $ G $, by $ E $ the exceptional divisor and by $ \overline{D_{100}}$, $ \overline{D_{101}}$, $ \overline{D_{110}}$, $ \overline{D_{111}}$ the strict transforms of $ D_{100}$, $ D_{101} $, $ D_{110}  $, $D_{111} $, respectively.

	Let $ \xymatrix{\overline{f}: \overline{X} \ar[r] & \overline{Y}}  $ be the $ \mathbb{Z}^3_2- $cover with  branch locus
	\begin{align*}
		B =\overline{D_{100}} + \overline{D_{101}} + \overline{D_{110}} + \overline{D_{111}}.
	\end{align*}
	The building data is as follows:
	$$	
	\begin{tabular}{l r r r}
		$ L_{100} \equiv$ & $ \left( m+2 \right)\overline{F}   $ & $ +\left( n+2\right)\overline{G} $& $ -2E $\\
		$ L_{010} \equiv$ & $m\overline{F}   $ &$ +n\overline{G}$& $ -E $\\
		$ L_{110} \equiv$ & $2\overline{F}   $ &$ +2\overline{G}$& $ -E $\\
		$ L_{001} \equiv$ & $ \overline{F}   $ & $ +\left( n+1\right)\overline{G} $& $ -E $\\
		$ L_{101} \equiv$ & $ \left( m+1\right)\overline{F}   $ & $ +\overline{G} $& $ -E $\\
		$ L_{011} \equiv$ & $ \left( m+1\right)\overline{F}   $ & $ +\overline{G} $& $ -E $\\
		$ L_{111} \equiv$ & $ \overline{F}   $ & $ +\left( n+1\right)\overline{G} $& $ -E $.
	\end{tabular} 
	$$
	\noindent
	Then the surface $\overline{X}$ is smooth and satisfies
	\begin{align*}
	2K_{\overline{X}} &\equiv f^*\left( 2m\overline{F} + 2n\overline{G} -2E \right) \\
					  &\equiv f^*\left( \left( 2m-1\right) \overline{F} + \left( 2n-1\right) \overline{G} +\overline{F_0} + \overline{G_0} \right)
	\end{align*}
	\noindent
	where $ \overline{F_0} $, $ \overline{G_0} $ are the strict transforms of the fibers $ F $, $ G $ passing through $ P $, respectively. Because the canonical divisor $ K_{\overline{X}} $ is nef and big, the surface $\overline{X}$ is minimal and of general type. Moreover, the surface $\overline{X}$ has the invariants:  
	\begin{align*}
		K_{\overline{X}}^2 = 16mn-8, p_g\left(\overline{X} \right)  = 2mn+2, q\left(\overline{X} \right) = 0.
	\end{align*}

	As in Section \ref{The canonical map and the canonical image}, the canonical map $ \varphi_{\left| K_{\overline{X}} \right| } $ is the composition of the quotient map $ \xymatrix{\overline{X} \ar[r]& Z } $ with the canonical map $\varphi_{\left| K_{Z} \right| }$ of $ Z $, where $ Z:= \overline{X}/\left\langle 0,0,1\right\rangle $. Moreover, the following diagram commutes
	$$
	\xymatrix{\overline{X} \ar[0,3]^{\mathbb{Z}_2^3}_{\overline{f}} \ar[1,1] \ar[2,1]_{\varphi_{\left| K_{\overline{X}} \right| }}&&& Y\\
		&Z \ar[-1,2]_{\mathbb{Z}_2^2} \ar[1,0]^{\varphi_{\left| K_{Z} \right| }}&&\\
		&\mathbb{P}^{2mn+1}&&} 	
	$$
	\noindent
	The eigenspace decomposition of $ H^{0}\left( K_Z\right)  $ is as follows:
	\begin{align*}
	H^{0}\left( K_Z\right) = H^{0}\left( \mathcal{O}_{Y}\right) \oplus H^{0}\left( \left( m-2\right) \overline{F} + \left( n-2\right) \overline{G}\right) \oplus H^{0}\left( m\overline{F}+n\overline{G} -E\right).
	\end{align*}	
	\noindent
	 Similarly to Section \ref{The canonical map and the canonical image}, one checks that $ \image\left( \varphi_{\left| K_{\overline{X}} \right| }\right)  $ is a surface of general type of degree $ 8mn-4$ and that $ \varphi_{\left| K_{\overline{X}}\right| } $ is a morphism of degree 2. Therefore, we obtain the surfaces described  in the second row of Theorem \ref{the main theorem}.	
		
	\begin{Remark}
      {\rm Taking  the  $ \mathbb{Z}_2^3 $ cover  of  $ \mathbb{P}^1 \times \mathbb{P}^1 $ with the above branch locus, we would obtain a singular surface with a Gorenstein elliptic singularity whose minimal resolution is an elliptic curve with self-intersection $-8$, (cf. \cite[ \rm No. 4.4, Table 1, Section 3.3]{MR2956036}).  The surface $ \overline{X} $ is the minimal resolution of this singular surface and  the pullback of the exceptional divisor $ E  $ is exactly the elliptic curve.}
	\end{Remark}
	\subsubsection{Variation 2}
	In this second variation of the construction given in Section \ref{the original construction of surface with d = 2} we again impose a quadruple point to the branch locus but in a different way.\\
	
	We take again a branch locus as before: $ D_{100}\in \left|  2F +2G\right|  $,  $ D_{101} $ $\in \left|  2F +2G\right|  $, $ D_{110} \in \left|  2mF \right|  $ (with $m\geq 2$) and $ D_{111}$ $\in \left|  2nG \right|  $ (with $n\geq 2$) , but now we require that $ D_{100} $ has  as unique singular point a node $ P $  and that   $ D_{101} $   and $ D_{110}$ go through the point $P$ in such a way that  the total branch $B$ has an ordinary  quadruple point.  As before we require that   $ D_{101} $, $ D_{110}$ and $D_{100}$ are smooth, that all intersections are transversal and that $P$ is the only point where more than two of the $D_\sigma$ meet.  
	 		
	Denote again  by $ \overline{Y} $ the blow up of $ \mathbb{P}^1 \times \mathbb{P}^1 $ at $ P $, by $ \overline{F} $, $ \overline{G} $ the pullback of general fibers $ F $, $ G $, by $ E $ the exceptional divisor and by $ \overline{D_{100}}$, $ \overline{D_{101}}$, $ \overline{D_{110}}$, $ \overline{D_{111}}$ the strict transforms of $ D_{100}$, $ D_{101} $, $ D_{110}  $, $D_{111} $, respectively.	
	
	Let $ \xymatrix{\overline{f}: \overline{X} \ar[r] & \overline{Y}}  $ be the $ \mathbb{Z}^3_2- $cover with  branch locus
	\begin{align*}
		B =\overline{D_{011}} +\overline{D_{100}} + \overline{D_{101}} + \overline{D_{110}} + \overline{D_{111}},
	\end{align*}
	\noindent
	where $ \overline{D_{011}}=E $. The building data is as follows:
	$$	
	\begin{tabular}{l r r r}
	$ L_{100} \equiv$ & $ \left( m+2 \right)\overline{F}   $ & $ +\left( n+2\right)\overline{G} $& $ -2E $\\
	$ L_{010} \equiv$ & $m\overline{F}   $ &$ +n\overline{G}$& $  $\\
	$ L_{110} \equiv$ & $2\overline{F}   $ &$ +2\overline{G}$& $ -E $\\
	$ L_{001} \equiv$ & $ \overline{F}   $ & $ +\left( n+1\right)\overline{G} $& $  $\\
	$ L_{101} \equiv$ & $ \left( m+1\right)\overline{F}   $ & $ +\overline{G} $& $ -E $\\
	$ L_{011} \equiv$ & $ \left( m+1\right)\overline{F}   $ & $ +\overline{G} $& $ -E $\\
	$ L_{111} \equiv$ & $ \overline{F}   $ & $ +\left( n+1\right)\overline{G} $& $ -E $.
	\end{tabular} 
	$$
	\noindent
	\noindent
	Then the surface $\overline{X}$ is smooth and satisfies
	\begin{align*}
	2K_{\overline{X}} &\equiv f^*\left( 2m\overline{F} + 2n\overline{G} -E \right) \\
					  &\equiv f^*\left( \left( 2m-1\right) \overline{F} + 2n\overline{G} +\overline{F_0} \right), 
	\end{align*}
	\noindent
	where $ \overline{F_0} $ is the strict transforms of the fiber $ F $ passing through $ P $. Because the canonical divisor $ K_{\overline{X}} $ is nef and big, the surface $\overline{X}$ is minimal and of general type. Moreover, the surface $\overline{X}$ has the invariants:
	\begin{align*}
		K_{\overline{X}}^2 = 16mn-2, p_g\left( \overline{X}\right) = 2mn+2, q\left( \overline{X}\right) = 0.
	\end{align*}
	
	As in Section \ref{The canonical map and the canonical image}, the canonical map $ \varphi_{\left| K_{\overline{X}} \right| } $ is the composition of the quotient map	$\xymatrix{\overline{X} \ar[r]^{h}& Z } $ with the canonical map $\varphi_{\left| K_{Z} \right| }$ of $ Z $, where $ Z:= \overline{X}/\left\langle 0,0,1\right\rangle $. Furthermore, the following diagram commutes
	$$
	\xymatrix{\overline{X} \ar[0,3]^{\mathbb{Z}_2^3}_{\overline{f}} \ar[1,1]^{h} \ar@{.>}[2,1]_{\varphi_{\left| K_{\overline{X}} \right| }}&&& Y\\
		&Z \ar[-1,2]_{\mathbb{Z}_2^2} \ar@{.>}[1,0]^{\varphi_{\left| K_{Z} \right| }}&&\\
		&\mathbb{P}^{2mn+1}&&} 	
	$$
	\noindent
	The eigenspace decomposition of $ H^{0}\left( K_Z\right)  $ is as follows:
	\begin{align*}
	H^{0}\left( K_Z\right) = H^{0}\left( \mathcal{O}_{Y}\right) \oplus H^{0}\left( \left( m-2\right) \overline{F} + \left( n-2\right) \overline{G} +E\right) \oplus H^{0}\left( m\overline{F}+n\overline{G} -E\right).
	\end{align*}	
	\noindent
	Similarly to Section \ref{The canonical map and the canonical image}, one checks that $ \image\left( \varphi_{\left| K_{\overline{X}} \right| }\right)  $ is a surface of general type and that $ \varphi_{\left| K_{\overline{X}}\right| } $ is a map of degree $ 2 $.\\
	
	Now we show that the canonical linear system $ \left| K_{\overline{X}}\right| $ has base points. First, we claim that 
	\begin{align}\label{Condition 1 of fixed point}
	K_{Z} \equiv g^{*}\left( m\overline{F}+n\overline{G} -E\right) + \tilde{E},
	\end{align}
	where the map $ \xymatrix{g :Z \ar[r] & Y} $ is the $ \mathbb{Z}_2^2- $cover with building data $ \left\lbrace L_{100}, L_{010}, L_{110} \right\rbrace $ and the divisor $ \tilde{E}: = \frac{1}{2}g^{*}\left( E \right) $ is an elliptic curve with self-intersection $ \tilde{E}^2=-1 $. In fact, the intermediate surface $ Z $ can be obtained by taking the iterated double covers
	\begin{align*}
	\xymatrix{Z \ar[r]^{g_2} & Z_1 \ar[r]^{g_1} & Y}
	\end{align*}
	\noindent
	where $ g_1 $ is the double cover branched on $ \overline{D_{100}} + \overline{D_{101}} + \overline{D_{110}} + \overline{D_{111}}$ and $ g_2 $ is the double cover branched on the pullback $ g_1^{*}\left( E\right)  $ and the nodes coming from the intersections between $ \overline{D_{100}} + \overline{D_{101}} $ and $ \overline{D_{110}} + \overline{D_{111}}$. So, we get
	\begin{align*}
	K_{Z_1} \equiv g_1^{*}\left( m\overline{F}+n\overline{G} -E\right).
	\end{align*}
	\noindent
	Because the elliptic curve $g_1^{*}\left( E \right) $ and the nodes are disjoint, we obtain the formula (\ref{Condition 1 of fixed point}).

	On the other hand, by the projection formula, we get
	\begin{align*}
	h^{0}\left( Z, g^{*}\left(m\overline{F}+n\overline{G} -E \right) \right) = 2mn+1
	\end{align*}
	\noindent
	This implies that
	\begin{align}\label{Condition 2 of fixed point}
	h^{0}\left( Z, K_Z \right) \ne h^{0}\left( Z, g^{*}\left(m\overline{F}+n\overline{G} -E \right) \right)
	\end{align}
	
	By $ (\ref{Condition 1 of fixed point}) $ and $ (\ref{Condition 2 of fixed point}) $, the elliptic curve $ \tilde{E} $ is not in the fixed part of $ \left|  K_{Z} \right|  $. On the other hand, because the divisor $ \left. K_{Z} \right|_{\tilde{E}} $ on the elliptic curve $ \tilde{E} $ has degree $ 1 $, the linear system $ \left| \left. K_{Z} \right|_{\tilde{E}}  \right|  $ has a simple base point. Therefore, the linear system $ \left| K_{Z}  \right|  $ has a simple base point.
	
	The double cover $ \xymatrix{\overline{X} \ar[r]& Z} $ ramifies on $ 8mn +12 $ nodes coming from the intersection points between $ \overline{D_{100}} $ and $ \overline{D_{101}} $, and $ \overline{D_{110}} $ and $ \overline{D_{111}} $. Because
	\begin{align*}
	K_{\overline{X}} \equiv h^{*}\left( K_{Z} \right) 
	\end{align*}
	\noindent
	and $ h^{0}\left( \overline{X}, K_{\overline{X}} \right)  = h^{0}\left( Z, K_{Z} \right)  $, the linear system $ \left| K_{\overline{X}}  \right|  $ possesses two simple base points on the pullback of $ \tilde{E} $. Moreover, the image $ \image\left( \varphi_{\left| K_{\overline{X}} \right| }\right)$ is a surface of degree $ 8mn-2 $.

	\begin{Remark}
	
	{\rm Taking the  $ \mathbb{Z}_2^3 $ cover  of  $ \mathbb{P}^1 \times \mathbb{P}^1 $ with the above branch locus, we would obtain a singular surface with a Gorenstein elliptic singularity whose minimal resolution is an elliptic curve with self-intersection $-2$. (cf. \cite[ \rm No. 4.8, Table 1, Section 3.3]{MR2956036}).  The surface $ \overline{X} $ is the minimal resolution of this singular surface and  the inverse image of the exceptional divisor $ E  $ is exactly the elliptic curve.}
	
	\end{Remark}

	\begin{Remark}
	In the three above constructions, if $ n = 1 $ or $ m=1 $, the canonical map $ \varphi_{\left| K_X\right| }  $ is a $ 4:1 $ map onto a rational ruled surface. 
	\end{Remark}

	\noindent
	In fact, if $ n = 1 $ or $ m=1 $, then the quotient surface $ Z/\left\langle 0,1,0\right\rangle  $ has geometric genus $ 0 $. This implies that the canonical map $ \varphi_{\left| K_{Z} \right| }$ has degree $ 2 $. Therefore, $ \degree\left( \varphi_{\left| K_X\right|} \right) = 4 $.\\

\begin{Acknowledgments}
The author is deeply indebted to Margarida Mendes Lopes for all her help. The author is supported by Funda\c{c}\~{a}o para a Ci\^{e}ncia e Tecnologia (FCT), Portugal under the framework of the program Lisbon Mathematics PhD (LisMath), Programa de Doutoramento FCT - PD/BD/113632/2015.
\end{Acknowledgments}	


\Addresses
\end{document}